# Traveling Wave Solutions for Space-Time Fractional Nonlinear Evolution Equations

M. G. Hafez[1] and Dianchen Lu[2*]

[1]Department of Mathematics, Chittagong University of Engineering and Technology, Chittagong-4349, Bangladesh. E-mail: golam_hafez@yahoo.com

[2]Faculty of Science, Jiangsu University, Zhenjiang, Jiangsu-212013, China.

*E-mail: dclu@ujs.edu.cn

## Abstract

Space-time fractional nonlinear evolution equations have been widely applied for describing various types of physical mechanism of natural phenomena in mathematical physics and engineering. The proposed generalized exp (-Φ (ξ))-expansion method along with the Jumarie's modified Riemann-Liouville derivatives is employed to carry out the integration of these equations, particularly space-time fractional Burgers equations, space-time fractional foam drainage equation and time fractional fifth order Sawada-Kotera equation. The traveling wave solutions of these equations are appeared in terms of the hyperbolic, trigonometric, exponential and rational functions. It has been shown that the proposed technique is a very effectual and easily applicable in investigation of exact traveling wave solutions to the fractional nonlinear evolution equations arises in mathematical physics and engineering.

**Keywords:** Space-time fractional nonlinear evolution equations; Traveling wave solutions; Modified Riemann-Liouville derivatives; Generalized exp (-Φ (ξ))-expansion method.

## 1. Introduction

The world around us is actually nonlinear and hence nonlinear evolution equations (NLEEs) are widely used as models and its versatile application in various fields of natural sciences [1, 2]. Nonlinear fractional partial differential equations (FPDEs) are a special class of NLEEs that have been focused several studies due to their frequent appearance in many application in physics, chemistry, biology, polymeric materials, electromagnetic, acoustics, neutron point kinetic model, vibration and control, signal and image processing, fluid dynamics and so on [3-7]. Due to its potential applications, researchers have devoted considerable effort to study the explicit and numerical solutions of nonlinear FPDEs. In order to understand the nonlinear physical mechanism of natural phenomena and further application in practical life, it is important to find more exact traveling wave solutions to the NLEEs. Already a large number of methods have applied to seek traveling wave solutions to nonlinear FPDEs, such as the fractional first integral method [8, 9], the fractional sub-equation method [10,11], the (*G′/G*)-expansion method [12-14], the improved (*G′/G*)-expansion method [15], the functional variable method [16], the fractional modified trial equation method [17], the extended spectral method [18], the variational iteration method[19] and so on. Li and He [20, 21] have proposed a fractional complex transformation to convert fractional differential equations into ordinary differential equations (ODEs). As a result, all analytical methods devoted to advance calculus can be easily applied to the fractional differential equations.

Recently, many authors [22-24] have applied the exp(-Φ(ξ))-expansion method to find the traveling wave solutions of the nonlinear PDEs arises in various fields as mention earlier. It is shown that the exp(-Φ(ξ))-expansion method according to the nonlinear ordinary differential equation (ODE) $\Phi'(\xi) = \exp(-\Phi(\xi)) + \mu \exp(\Phi(\xi)) + \lambda, \lambda, \mu \in \Re$ have been given few comprehensive solutions to the nonlinear PDEs. Very recently, Hafez and Akbar [25] have applied the exp(-Φ(ξ))-expansion method to solve strain wave equation appeared in microstructure solids by considering $\Phi'(\xi) = \exp(-\Phi(\xi)) + \mu \exp(\Phi(\xi)) + \lambda$, $\Phi'(\xi) = -\sqrt{\lambda + \mu \exp(-\Phi(\xi))^2}$ and $\Phi'(\xi) = -\lambda \exp(-\Phi(\xi)) - \mu \exp(\Phi(\xi))$ as auxiliary differential equations. In order to obtain the standard form of this method, we have used the nonlinear differential equation $\Phi'(\xi) = p\exp(-\Phi(\xi)) + q\exp(\Phi(\xi)) + r$ as auxiliary equation for finding more comprehensive solutions to nonlinear PDEs, so-called generalized exp(-Φ(ξ))-expansion method. Therefore, the purpose of this paper is to present the proposed generalized exp (-Φ (ξ))-expansion method and apply this method to construct the exact traveling wave solutions of the space-time fractional Burgers equation [26], space-time fractional coupled Burger's equation [27],space-time fractional foam drainage equation and time fractional fifth order Sawada-Kotera (SK) equation. The proposed method according to the auxiliary nonlinear ODE provides much more comprehensive results and easily applicable to solve the NLEEs. Moreover, we have employed this method for finding more comprehensive exact traveling wave solutions to the nonlinear FPDEs. Sometimes this method can be given solutions in disguised versions of known solutions that may be obtained by other methods. The advantage of this method over the existing method is that it provides some new exact traveling wave solutions together with additional free parameters. Apart from the physical significance, the close-form solutions of NLEEs may be helpful the numerical solvers to compare the correctness of their results and help them in the stability analysis. The algebraic manipulation of this method with the help of algebraic software such as, Maple is much easier than the other accessible method.

The rest of the paper has been prepared as follows: In section 2, the proposed generalized $\exp(-\Phi(\xi))$-expansion method is discussed in details. The section 3 presents the application of this method to construct the exact traveling wave solutions of the nonlinear FPDEs. The advantages of the proposed method and comparison with others methods is given in section 4. Conclusions have been drawn in Section 5.

## 2. Representation of the method

This section presents the brief descriptions of the proposed method for finding the exact traveling to the nonlinear FPDEs.

One can consider the nonlinear FPDEs to the following form:

$$f(u, D_t^{\alpha}u, D_x^{\alpha}u, D_x^{\beta}u, D_t^{\alpha} D_t^{\beta}u, D_x^{\alpha} D_x^{\beta}u, \ldots\ldots) = 0, \alpha > 0, \beta < 1, \tag{1}$$

where, $f$ is a function of $u(x,t)$ and its partial fractional derivatives, in which higher order derivatives and nonlinear terms are involved.

Let us consider the traveling wave variable transformation as

$$u(x,t) = U(\xi),\ \xi = \frac{kx^{\beta}}{\Gamma(1+\beta)} \pm \frac{ct^{\alpha}}{\Gamma(1+\alpha)}, \tag{2}$$

where $k$ and $c$ are constants.

The Jumarie's modified Riemann-Liouville [28] derivatives of order $\alpha$ can be defined by the following expression

$$D_z^{\alpha} f(z) = \begin{cases} \dfrac{1}{\Gamma(1-\alpha)} \dfrac{d}{dz} \displaystyle\int_0^z (z-\xi)^{-\alpha} (f(\xi) - f(0)) d\xi;\ 0 < \alpha < 1 \\ (f^{(n)}(z))^{\alpha - n};\ n \le \alpha < n+1, n \ge 1. \end{cases} \tag{3}$$

Moreover, one can define the modified Riemann-Liouville derivative [29, 30] as

$$D_z^{\alpha} z^{\gamma} = \frac{\Gamma(1+\gamma)}{\Gamma(1+\gamma-\alpha)} z^{\gamma-\alpha}, \gamma > 0, \tag{4}$$

Using Eq. (2) and (4), Eq. (1) can be converted as a nonlinear ordinary differential equation (ODE) for $U = U(\xi)$:

$$F(U, U', U'', U''', \cdots\cdots\cdots) = 0, \tag{5}$$

where, $F$ is a function of $U, U', U'', U''', \cdots\cdots\cdots$ and its derivatives point out the ordinary derivatives with respect to $\xi$.

One can be considered the traveling wave solutions of Eq. (5) as

$$u(\xi)=\sum_{i=0}^{N}A_i\left(e^{-\Phi(\xi)}\right)^i,\ A_N\neq 0 \tag{6}$$

where the coefficients $A_i\,(0\le i\le N)$ are constants to be evaluated and $\Phi=\Phi(\xi)$ satisfies the following first order nonlinear ordinary differential equation:

$$\Phi'(\xi)=p\,e^{-\Phi(\xi)}+q\,e^{\Phi(\xi)}+r, \tag{7}$$

Balancing the higher order derivative with the nonlinear terms of the highest order that appeared in Eq. (5), one can get the value of the positive integer $N$. On the other hand, if the degree of $U=U(\xi)$ is $D[U(\xi)]=n$, then the degree of the other expressions can be found by the following formulae:

$$D[\frac{d^N U(\xi)}{d\xi^N}]=N+p,\quad D[u^N\left(\frac{d^K U(\xi)}{d\xi^K}\right)^S]=nN+S(n+k). \tag{8}$$

Substituting Eq. (6) into Eq. (5) and using (7) rapidly, one can obtain a system of algebraic equations for $A_i\,(0\le i\le N),k,p,q,r,c$. With the help of symbolic computation, such as Maple, one can evaluate the obtaining system and find out the values $A_i\,(0\le i\le N),k,p,q,r,c$. It is notable that equation (7) has the following general solutions:

**Type 1**: when $p=1,$ one obtains

$$\Phi(\xi)=\ln\left(-\frac{\sqrt{r^2-4q}\tanh\left\{0.5\sqrt{r^2-4q}(\xi+\xi_0)\right\}}{2q}-\frac{\gamma}{2q}\right),q\neq 0,r^2-4q>0, \tag{9a}$$

$$\Phi(\xi)=\ln\left(\frac{\sqrt{r^2-4q}\tan\left\{0.5\sqrt{r^2-4q}(\xi+\xi_0)\right\}}{2q}-\frac{\gamma}{2q}\right),\ q\neq 0,r^2-4q<0, \tag{9b}$$

$$\Phi(\xi)=-\ln\left(\frac{r}{\exp(r(\xi+\xi_0))-1}\right),q=0,r\neq 0,\ r^2-4q>0, \tag{9c}$$

$$\Phi(\xi)=\ln\left(-\frac{2(r(\xi+\xi_0)+2)}{r^2(\xi+\xi_0)}\right),q\neq 0,r\neq 0,\ r^2-4q=0, \tag{9d}$$

**Type 2**: when $r=0,$ one obtains

$$\Phi(\xi)=\ln\left(\sqrt{\frac{p}{q}}\tan\left[\sqrt{pq}\ (\xi+\xi_0)\right]\right),p>0,q>0, \tag{9e}$$

$$\Phi(\xi) = \ln\left( -\sqrt{\frac{p}{q}} \cot\left[\sqrt{pq}\ (\xi + \xi_0)\right]\right),\ p > 0, q > 0, \tag{9f}$$

$$\Phi(\xi) = \operatorname{sgn}(p) \ln\left( -\sqrt{\frac{-p}{q}} \tanh\left[\sqrt{-pq}\ (\xi + \xi_0)\right]\right),\ pq < 0, \tag{9g}$$

$$\Phi(\xi) = \operatorname{sgn}(p) \ln\left( -\sqrt{\frac{-p}{q}} \coth\left[\sqrt{-pq}\ (\xi + \xi_0)\right]\right),\ pq < 0\,, \tag{9h}$$

**Type 3:** when $q = 0$ and $r = 0$, one obtains

$$\Phi(\xi) = \ln\{p(\xi + \xi_0)\}, \tag{9i}$$

where $\xi_0$ is the integrating constant.

Finally, we are obtained the multiple explicit solutions of nonlinear FPDE (1) by combining the equations (2), (6) and (9).

## 3. Applications of the method

This section presents four examples to illustrate the applicability of the proposed method to solve the space-time FDEs.

### 3.1 Space-time fractional Burgers equation

Let us consider the space-time Burgers equation as follows

$$\frac{\partial^{\alpha} u}{\partial t^{\alpha}} - 2u\frac{\partial^{\beta} u}{\partial x^{\beta}} - A\frac{\partial^{2\beta} u}{\partial x^{2\beta}} = 0,\ \alpha > 0,\ \beta < 1,\ \ x > 0,\ \ t > 0,\ \ A \in \Re \tag{10}$$

When $\alpha = \beta = 1$, the fractional equation (10) can be reduced to the well known Burgers equation. The space-time FDE (10) have appeared a mathematical model equation not only in fluid flow [31] but also it can be applied for describing various types of physical phenomena in the field of gas dynamics, heat conduction, elasticity, continuous stochastic processes, etc.

Introducing the following transformation in eq. (10),

$$u(x,t) = w(\xi),\ \ \xi = \frac{kx^{\beta}}{\Gamma(1+\beta)} - \frac{ct^{\alpha}}{\Gamma(1+\alpha)}, \tag{11}$$

where $k$ and $c$ are constants, one can rewrite the equation (11) to the following nonlinear ODE as follows:

$$cw' + 2kww' + Ak^2 w'' = 0, \tag{12}$$

where primes denote the differentiation with respect to $\xi$. Integrating the eq. (12) once and setting the constant of integration to zero for time homogeneity, we get

$$cw + kw^2 + Ak^2 w' = 0, \tag{13}$$

Using the balancing principle between $w'$ and $w^2$ in eq. (13) gives $N = 1$. Therefore the solution of (13) can be written as

$$w(\xi) = A_0 + A_1 e^{-\Phi(\xi)} \tag{14}$$

where $A_0$, $A_1$ are constants to be determined later and $\Phi(\xi)$ satisfies the auxiliary nonlinear ODE (7).

Substituting eq. (14) into eq. (13) and using (7) frequently, the left-hand side of eq. (13) becomes a polynomial in $e^{-\Phi(\xi)}$. Setting the coefficients of this polynomial to zero yields a system of algebraic equations as follows:

$$\left\{cA_0 - Ak^2 A_1 q + kA_0^2 = 0, cA_1 + 2kA_0 A_1 - Ak^2 A_1 r = 0,\ kA_1^2 - Ak^2 A_1 p = 0\right\} \tag{15}$$

Solving the resulting algebraic equations (15), we have

$$\left\{A_0 = Ak\left(\frac{r \pm \sqrt{r^2 - 4pq}}{2}\right), A_1 = pAk, k = k, c = -Ak^2\left(r \pm \sqrt{r^2 - 4pq}\right)\right\} \tag{16}$$

where $k$, $p, q$ and $r$ are arbitrary constants.

Combing the solutions of Eq. (7), (11), (14) and (16), one can obtain the following explicit solutions to the space-time fractional Burgers equation (10):

**For Type 1:**

$$u_{1_1}(\xi) = Ak\left\{\frac{r \pm \sqrt{r^2 - 4q}}{2} - \left(\frac{2q}{\sqrt{r^2 - 4q}\tanh(0.5\sqrt{r^2 - 4q}(\xi + \xi_0)) + r}\right)\right\}, \quad q \neq 0, r^2 - 4q > 0, \tag{17}$$

$$u_{1_2}(\xi) = Ak\left\{\frac{r \pm \sqrt{r^2 - 4q}}{2} + \left(\frac{2q}{\sqrt{-(r^2 - 4q)}\tan(0.5\sqrt{-(r^2 - 4q)}(\xi + \xi_0)) - r}\right)\right\}, \quad q \neq 0, r^2 - 4q < 0, \tag{18}$$

$$u_{1_3}(\xi) = Ak\left\{\frac{r \pm \sqrt{r^2 - 4q}}{2} + \left(\frac{r}{\exp(r(\xi+\xi_0)) - 1}\right)\right\}, q = 0,\ r^2 - 4q > 0 \qquad (19)$$

$$u_{1_4}(x,t) = Ak\left\{\frac{r \pm \sqrt{r^2 - 4q}}{2} - \left(\frac{r^2(\xi+\xi_0)}{2r(\xi+\xi_0)) + 4}\right)\right\}, q \neq 0, r \neq 0,\ r^2 - 4q = 0, \qquad (20)$$

where, $\xi = \dfrac{kx^{\beta}}{\Gamma(1+\beta)} + \dfrac{Ak^2(r \pm \sqrt{r^2 - 4q})}{\Gamma(1+\alpha)} t^{\alpha}$.

**For Type 2:**

$$u_{1_5}(x,t) = Ak\sqrt{-pq}\left\{\pm 1 - \tanh\left[\sqrt{-pq}\left(\frac{kx^{\beta}}{\Gamma(1+\beta)} - \frac{\mp 2k^2 A\sqrt{-pq}t^{\alpha}}{\Gamma(1+\alpha)} + \xi_0\right)\right]\right\}, p < 0, p > 0, \qquad (21)$$

$$u_{1_6}(x,t) = Ak\sqrt{-pq}\left\{\pm 1 - \coth\left[\sqrt{-pq}\left(\frac{kx^{\beta}}{\Gamma(1+\beta)} - \frac{\mp 2k^2 A\sqrt{-pq}t^{\alpha}}{\Gamma(1+\alpha)} + \xi_0\right)\right]\right\}, p < 0, p > 0,, \qquad (22)$$

$$u_{1_7}(x,t) = Ak\left\{\pm\sqrt{-pq} - \sqrt{pq}\tan\left[\sqrt{pq}\left(\frac{kx^{\beta}}{\Gamma(1+\beta)} - \frac{\mp 2k^2 a\sqrt{-pq}t^{\alpha}}{\Gamma(1+\alpha)} + \xi_0\right)\right]\right\}, p > 0, q > 0, \qquad (23)$$

$$u_{1_8}(x,t) = Ak\sqrt{pq}\left\{\pm 1 + \cot\left[\sqrt{pq}\left(\frac{kx^{\beta}}{\Gamma(1+\beta)} - \frac{\mp 2k^2 a\sqrt{pq}t^{\alpha}}{\Gamma(1+\alpha)} + \xi_0\right)\right]\right\}, pq > 0, \qquad (24)$$

**For Type 3:**

$$u_{1_9}(x,t) = Ak\left\{\frac{\gamma \pm \sqrt{\gamma^2 - 4\beta}}{2} + \left(\frac{kx^{\beta}}{\Gamma(1+\beta)} + \frac{Ak^2(2\gamma \pm \sqrt{\gamma^2 - 4\beta})}{\Gamma(1+\alpha)} t^{\alpha}\xi + \xi_0\right)^{-1}\right\}, q = 0,\ r = 0, \qquad (25)$$

## 3.2 Space-time fractional coupled Burgers equation

Let us consider the space-time fractional coupled Burgers equation as follows

$$\left.\begin{aligned} \frac{\partial^{\alpha} u}{\partial t^{\alpha}} - \frac{\partial^{2\alpha} u}{\partial x^{2\alpha}} + 2u\frac{\partial^{\alpha} u}{\partial x^{\alpha}} + L\frac{\partial^{\alpha}(uv)}{\partial x^{\alpha}} = 0 \\ \frac{\partial^{\alpha} v}{\partial t^{\alpha}} - \frac{\partial^{2\alpha} v}{\partial x^{2\alpha}} + 2v\frac{\partial^{\alpha} v}{\partial x^{\alpha}} + M\frac{\partial^{\alpha}(uv)}{\partial x^{\alpha}} = 0 \end{aligned}\right\}. \qquad (26)$$

The coupled fractional equations have appeared as model equation in mathematical physics, which is derived by Esipov [32]. It is very significant that the system is a simple model of sedimentation or evolution of scaled volume concentrations of two kinds of particles in fluid suspensions or colloids, under the effect of gravity [33]. The constants $L$ and $M$ depend on the system parameters such as the Peclet number, the Stokes velocity of particles due to gravity, and the Brownian diffusivity.

One can introduce the following transformation

$$u(x,t)=u(\xi),v(x,t)=v(\xi),\ \xi=\frac{x^{\alpha}}{\Gamma(1+\alpha)}+\frac{ct^{\alpha}}{\Gamma(1+\alpha)} \tag{27}$$

where $c$ is a constant, Eq. (26) can be converted to the following form

$$\left.\begin{aligned} &cu'-u''+2uu'+L(uv)'=0\\ &cv'-v''+2vv'+M(uv)'=0\end{aligned}\right\} \tag{28}$$

where primes denote the differentiation with respect to $\xi$.

According to the balancing principle, the solution of the system of eq. (28) can be expressed by a polynomial in $e^{-\Phi(\xi)}$ as follows:

$$\left.\begin{aligned} &u(\xi)=A_0+A_1e^{-\Phi(\xi)}\\ &v(\xi)=B_0+B_1e^{-\Phi(\xi)}\end{aligned}\right\}, \tag{29}$$

where $A_0,\ A_1, B_0,\ B_1$ are constants to be determined later and $\Phi(\xi)$ satisfies the equation (7).

By substituting eq. (29) into eq. (28) and using (7) frequently, one can obtain the following system of algebraic equations by setting the coefficients of the polynomial in $e^{-\Phi(\xi)}$ to zero:

$$\left.\begin{aligned} &-A_1qr-2A_1A_0q-LB_1A_0q-LA_1B_0q-cA_1q=0,\\ &-2A_1^2q-2A_1pq-A_1r^2-2A_1A_0r-cA_1r-LA_1B_0r-LA_0B_1r-2LA_1B_1q=0\\ &-2LB_1A_1r-2A_1A_0p-2A_1^2r-cA_1p-3A_1pr-LpB_1A_0-LpB_0A_1=0,-2A_1p^2-2A_1^2p-2LB_1A_1p=0,\\ &-B_1qr-2B_1B_0q-MB_1A_0q-MA_1B_0q-cB_1q=0,\\ &-2B_1^2q-2B_1pq-B_1r^2-2B_1B_0r-cB_1r-MB_1A_0r-MA_1B_0r-2MA_1B_1q=0\\ &-2MB_1A_1r-2B_1B_0p-2B_1^2r-cB_1p-3B_1pr-MpB_1A_0-MpB_0A_1=0,-2B_1p^2-2B_1^2p-2MB_1A_1p=0\end{aligned}\right\}. \tag{30}$$

Solving the resulting algebraic equations (30) with aid of symbolic computation, such as Maple, one obtains

$$\left\{c=-\frac{2LMB_0-r+Mr-2B_0}{-1+M},\ A_0=\frac{(-1+L)B_0}{(-1+M)},A_1=-\frac{p(-1+L)}{-1+LM},B_0=B_0,B_1=-\frac{p(-1+M)}{-1+LM}\right\} \tag{31}$$

where $B_0$, $p, q$ and $r$ are arbitrary constants.

By combing the equations (9), (27), (29) and (31), the space-time fractional coupled Burger's equation (26) has the following traveling wave solutions:

**For Type 1:**

$$\left.\begin{aligned} u_1(x,t) &= \frac{(-1+L)B_0}{(-1+M)} + \frac{(-1+L)}{-1+LM}\left(\frac{2q}{\sqrt{r^2-4q}\tanh(0.5\sqrt{r^2-4q}(\xi+\xi_0))+r}\right) \\ v_1(x,t) &= B_0 + \frac{(-1+M)}{-1+LM}\left(\frac{2\mu}{\sqrt{r^2-4q}\tanh(0.5\sqrt{r^2-4q}(\xi+\xi_0))+r}\right) \end{aligned}\right\}, q \neq 0, r^2-4q>0, \quad (32)$$

$$\left.\begin{aligned} u_2(x,t) &= \frac{(-1+L)B_0}{(-1+M)} - \frac{(-1+L)}{-1+LM}\left(\frac{2q}{\sqrt{4q-r^2}\tan(0.5\sqrt{4q-r^2}(\xi+\xi_0))-r}\right) \\ v_2(x,t) &= B_0 - \frac{(-1+M)}{-1+LM}\left(\frac{2q}{\sqrt{4q-r^2}\tan(0.5\sqrt{4q-r^2}(\xi+\xi_0))-r}\right) \end{aligned}\right\}, q \neq 0, r^2-4q<0, \quad (33)$$

$$\left.\begin{aligned} u_3(x,t) &= \frac{(-1+L)B_0}{(-1+M)} - \frac{(-1+L)}{-1+LM}\left(\frac{r}{\exp(r(\xi+\xi_0))-1}\right) \\ v_3(x,t) &= B_0 - \frac{(-1+M)}{-1+LM}\left(\frac{r}{\exp(r(\xi+\xi_0))-1}\right) \end{aligned}\right\}, q = 0,\ r^2-4q>0, \quad (34)$$

$$\left.\begin{aligned} u_4(x,t) &= \frac{(-1+L)B_0}{(-1+M)} + \frac{(-1+L)}{-1+LM}\left(\frac{r^2(\xi+\xi_0)}{2r(\xi+\xi_0))+4}\right) \\ v_4(x,t) &= B_0 + \frac{(-1+M)}{-1+LM}\left(\frac{r^2(\xi+\xi_0)}{2r(\xi+\xi_0))+4}\right) \end{aligned}\right\}, q \neq 0, r \neq 0,\ r^2-4q=0, \quad (35)$$

where, $\xi = \frac{x^\alpha}{\Gamma(1+\alpha)} + \frac{ct^\alpha}{\Gamma(1+\alpha)}$ and $c = -\frac{2LMB_0 - r + Mr - 2B_0}{-1+M}$.

**For Type 2:**

$$\left.\begin{aligned} u_5(x,t) &= \frac{(-1+L)B_0}{(-1+M)} + \frac{(-1+L)}{-1+LM}\sqrt{pq}\tan\left[\sqrt{pq}\left(\frac{x^\alpha}{\Gamma(1+\alpha)} + \frac{ct^\alpha}{\Gamma(1+\alpha)} + \xi_0\right)\right] \\ v_5(x,t) &= B_0 + \frac{(-1+M)}{-1+LM}\sqrt{pq}\tan\left[\sqrt{pq}\left(\frac{x^\alpha}{\Gamma(1+\alpha)} + \frac{ct^\alpha}{\Gamma(1+\alpha)} + \xi_0\right)\right] \end{aligned}\right\}, p>0, q>0 \quad (36)$$

$$\left.\begin{aligned} u_6(x,t) &= \frac{(-1+L)B_0}{(-1+M)} - \frac{(-1+L)}{-1+LM}\sqrt{pq}\cot\left[\sqrt{pq}\left(\frac{x^\alpha}{\Gamma(1+\alpha)} + \frac{ct^\alpha}{\Gamma(1+\alpha)} + \xi_0\right)\right] \\ v_6(x,t) &= B_0 - \frac{(-1+M)}{-1+LM}\sqrt{pq}\cot\left[\sqrt{pq}\left(\frac{x^\alpha}{\Gamma(1+\alpha)} + \frac{ct^\alpha}{\Gamma(1+\alpha)} + \xi_0\right)\right] \end{aligned}\right\}, p>0, q>0, \qquad (37)$$

$$\left.\begin{aligned} u_7(x,t) &= \frac{(-1+L)B_0}{(-1+M)} + \frac{(-1+L)}{-1+LM}\sqrt{-pq}\tanh\left[\sqrt{-\lambda\mu}\left(\frac{x^\alpha}{\Gamma(1+\alpha)} + \frac{ct^\alpha}{\Gamma(1+\alpha)} + \xi_0\right)\right] \\ v_7(x,t) &= B_0 + \frac{(-1+M)}{-1+LM}\sqrt{-pq}\tanh\left[\sqrt{-pq}\left(\frac{x^\alpha}{\Gamma(1+\alpha)} + \frac{ct^\alpha}{\Gamma(1+\alpha)} + \xi_0\right)\right] \end{aligned}\right\}, p<0, q>0, \qquad (38)$$

$$\left.\begin{aligned} u_8(x,t) &= \frac{(-1+L)B_0}{(-1+M)} + \frac{(-1+L)}{-1+LM}\sqrt{-pq}\coth\left[\sqrt{-pq}\left(\frac{x^\alpha}{\Gamma(1+\alpha)} + \frac{ct^\alpha}{\Gamma(1+\alpha)} + \xi_0\right)\right] \\ v_8(x,t) &= B_0 + \frac{(-1+M)}{-1+LM}\sqrt{-pq}\coth\left[\sqrt{-pq}\left(\frac{x^\alpha}{\Gamma(1+\alpha)} + \frac{ct^\alpha}{\Gamma(1+\alpha)} + \xi_0\right)\right] \end{aligned}\right\}, p<0, q>0, \qquad (39)$$

where $c = -\dfrac{2LMB_0 - 2B_0}{-1+M}$.

**For Type 3:**

$$\left.\begin{aligned} u_9(x,t) &= \frac{(-1+L)B_0}{(-1+M)} - \frac{(-1+L)}{-1+LM}\left(\frac{x^\alpha}{\Gamma(1+\alpha)} + \frac{ct^\alpha}{\Gamma(1+\alpha)} + \xi_0\right)^{-1} \\ v_9(x,t) &= B_0 - \frac{(-1+M)}{-1+LM}\left(\frac{x^\alpha}{\Gamma(1+\alpha)} + \frac{ct^\alpha}{\Gamma(1+\alpha)} + \xi_0\right)^{-1} \end{aligned}\right\}, \mu=0,\ \lambda=0, \qquad (40)$$

where, $c = -\dfrac{2LMB_0 - 2B_0}{-1+M}$.

### 3.3 Space-time fractional foam drainage equation

Let us consider the following space-time fractional foam drainage equation

$$\frac{\partial^\alpha V}{\partial t^\alpha} = \frac{1}{2}V\frac{\partial^{2\beta}V}{\partial t^{2\beta}} + 2V^2\frac{\partial^\beta V}{\partial x^\beta} + \left(\frac{\partial^\beta V}{\partial x^\beta}\right)^2, \ \alpha>0,\ \beta\le 1,\ t>0, \qquad (41)$$

It has been studied by many authors [15, 34]. This equation have appeared as a simple model for describing the flow of liquid through channels ( Plateau borders [35] ) and nodes ( intersection of four channels) between the bubbles, driven by gravity and capillarity [36].

Introducing the fractional complex transformation

$$V(x,t) = V(\xi),\ \ \xi = \frac{kx^{\beta}}{\Gamma(1+\beta)} + \frac{ct^{\alpha}}{\Gamma(1+\alpha)}, \tag{42}$$

where $k$ and $c$ are constants in eq. (41), Eq. (41) can be transformed as

$$-cV' + \frac{1}{2}k^2VV'' + 2kV^2V' + k^2V'^2 = 0, \tag{43}$$

According to the balancing principle, the solution of (43) can be written as

$$V(\xi) = A_0 + A_1e^{-\Phi(\xi)} \tag{44}$$

where $A_0$, $A_1$ are constants to be determined later and $\Phi(\xi)$ satisfies the nonlinear ODE (7).

By substituting eq. (44) into eq. (43) and using (7) commonly, the left-hand side of eq. (43) becomes a polynomial in $e^{-\Phi(\xi)}$. Collecting the coefficients of this polynomial to zero yields a system of algebraic equations in terms of $A_0, A_1, k, p, q, r, c$. The algebraic equations are overlooked for convenience. Solving the resulting algebraic equations with aid of symbolic computation, such as Maple, one get

$$\left\{A_0 = \frac{k}{2}r, A_1 = kp, k = k, c = -k^3pq + \frac{1}{4}k^3r^2\right\} \tag{45}$$

where $k$, $p, q$ and $r$ are arbitrary constants.

By combining the equations (9), (42), (44) and (45), the space-time fractional foam drainage equation (41) has the following traveling wave solutions:

**For Type 1:**

$$V_1(x,t) = \frac{k}{2}r - k\left(\frac{2q}{\sqrt{r^2-4q}\tanh(0.5\sqrt{r^2-4q}(\xi+\xi_0)) + r}\right),\ \ q \neq 0, r^2 - 4q > 0, \tag{46}$$

$$V_2(x,t) = \frac{k}{2}r + k\left(\frac{2q}{\sqrt{4q-r^2}\tan(0.5\sqrt{4q-r^2}(\xi+\xi_0)) - r}\right),\ \ q \neq 0, r^2 - 4q < 0, \tag{47}$$

$$V_3(x,t) = \frac{k}{2}r + k\left(\frac{r}{\exp(r(\xi+\xi_0)) - 1}\right), q = 0,\ r^2 - 4q > 0, \tag{48}$$

$$V_4(x,t) = \frac{k}{2}r - k\left(\frac{r^2(\xi+\xi_0)}{2r(\xi+\xi_0)) + 4}\right), q \neq 0, r \neq 0,\ r^2 - 4q = 0, \tag{49}$$

where, $\xi = \frac{kx^{\beta}}{\Gamma(1+\beta)} + \frac{-4k^3q + k^3r^2}{4\,\Gamma(1+\alpha)} t^{\alpha}$.

**For Type 2:**

$$V_5(x,t) = -k\sqrt{pq}\tan\left[\sqrt{pq}\left(\frac{kx^{\alpha}}{\Gamma(1+\beta)} + \frac{-k^3pqt^{\alpha}}{\Gamma(1+\alpha)} + \xi_0\right)\right], p>0, q>0, \quad (50)$$

$$V_6(x,t) = k\sqrt{pq}\cot\left[\sqrt{pq}\left(\frac{kx^{\alpha}}{\Gamma(1+\beta)} + \frac{-k^3pqt^{\alpha}}{\Gamma(1+\alpha)} + \xi_0\right)\right], p>0, q>0, \quad (51)$$

$$V_7(x,t) = -k\sqrt{-pq}\tanh\left[\sqrt{-pq}\left(\frac{kx^{\alpha}}{\Gamma(1+\beta)} + \frac{-k^3pqt^{\alpha}}{\Gamma(1+\alpha)} + \xi_0\right)\right], p<0, q>0, \quad (52)$$

$$V_8(x,t) = -k\sqrt{-pq}\coth\left[\sqrt{-pq}\left(\frac{kx^{\alpha}}{\Gamma(1+\beta)} + \frac{-k^3pqt^{\alpha}}{\Gamma(1+\alpha)} + \xi_0\right)\right], p<0, q>0, \quad (53)$$

**For Type 3:**

$$V_9(x,t) = k\left(\frac{kx^{\beta}}{\Gamma(1+\beta)} + \frac{-4k^3q + k^3r^2}{4\,\Gamma(1+\alpha)} t^{\alpha} + \xi_0\right)^{-1}, q=0, r=0, \quad (54)$$

### 3.4 Time fractional fifth order Sawada-Kotera (SK) equation

Finally, to illustrate more applicability of this proposed method, the following time fractional SK equation is considered:

$$\frac{\partial^{\alpha}U}{\partial t^{\alpha}} + 5U^2\frac{\partial U}{\partial x} + 5\frac{\partial U}{\partial x}\frac{\partial^2 U}{\partial x^2} + 5U\frac{\partial^3 U}{\partial x^3} + \frac{\partial^5 U}{\partial x^5} = 0,\ 0<\alpha\le 1. \quad (55)$$

The time fractional equation (55) appeared as a model equation in many physical instances that propagates in opposite directions. To obtain the analytical solutions of Eq. (55), one can introduce the following fractional complex transformation:

$$U(x,t) = w(\xi),\ \xi = k[x - \frac{ct^{\alpha}}{\Gamma(1+\alpha)}], \quad (56)$$

where $k$ and $c$ are constants to be evaluated later. According to the transformation (55), one can be obtained to the following nonlinear ODE:

$$-kc\frac{dw}{d\xi}+\frac{5}{3}k\frac{d}{d\xi}(w^3)+5k^3\frac{d}{d\xi}(w\frac{d^2w}{d\xi^2})+k^5\frac{d^5w}{d\xi^5}=0. \quad (57)$$

Integrating the Eq. (57) once and setting the constant of integration to zero for simplicity, one get

$$-cw+\frac{5}{3}w^3+5k^2w\frac{d^2w}{d\xi^2}+k^4\frac{d^4w}{d\xi^4}=0. \quad (58)$$

According to the balancing principle between $\frac{d^4w}{d\xi^4}$ and $w^3$ that are involved in eq. (58), one obtain $N=2.$ Therefore the solution of (58) according to the proposed method can be written as

$$w(\xi)=A_0+A_1e^{-\Phi(\xi)}+A_2e^{-2\Phi(\xi)} \quad (59)$$

where $A_0$, $A_1$ and $A_2$ are constants to be determined later. Substituting Eq. (59) into Eq. (58) and collecting the degree of the polynomial in $e^{-\Phi(\xi)}$ yields a system of algebraic nonlinear equations which are omitted for simplicity. Solving the resulting algebraic equations, one obtains

**Set 1:** $c=k^4(-8pqr^2+r^4+16p^2q^2), A_0=-6k^2pq, A_1=-6k^2pr, A_2=-6k^2p^2$ (60)

**Set2:** $c=(-\frac{5k^4r^2}{2}+10k^4pq)\chi-11k^4r^2qp-\frac{k^4r^4}{2}+52k^4p^2q^2, A_0=\chi k^2, A_1=-6k^2pr, A_2=-6k^2p^2$, (61)

where $\chi=\frac{-(60pq+15r^2)\pm\sqrt{105r^4-1680p^2q^2-840pqr^2}}{20}$, $p,q$ and $r$ are arbitrary constants.

Therefore, the time fractional SK equation according to the **Set 1** has the following traveling wave solutions:

**For Type 1:**

$$U_{1_1}(x,t)=-6k^2q+6k^2r\left(\frac{q}{\sqrt{r^2-4q}\tanh(0.5\sqrt{r^2-4q}(\xi+\xi_0))+r}\right)- \\ 6k^2\left(\frac{q}{\sqrt{r^2-4q}\tanh(0.5\sqrt{r^2-4q}(\xi+\xi_0))+r}\right)^2, r^2-4q>0, q\neq 0, \quad (62)$$

$$U_{1_2}(x,t) = -6k^2 q - 6k^2 r\left(\frac{q}{\sqrt{-(r^2-4q)}\tan(0.5\sqrt{-(r^2-4q)}(\xi+\xi_0)) - r}\right) - 6k^2\left(\frac{2q}{\sqrt{-(r^2-4q)}\tan(0.5\sqrt{-(r^2-4q)}(\xi+\xi_0)) - r}\right)^2, r^2-4q<0, \mu\neq 0, \quad (63)$$

$$U_{1_3}(x,t) = -6k^2\left(\frac{r^2}{e^{(r(\xi+\xi_0))}-1}\right) - 6k^2\left(\frac{r}{e^{(r(\xi+\xi_0))}-1}\right)^2, q=0, r\neq 0, (r^2-4q)>0 \quad (64)$$

$$U_{1_4}(x,t) = -6k^2 q + 6k^2\left(\frac{r^3(\xi+\xi_0)}{2r(\xi+\xi_0))+4}\right) - 6k^2\left(\frac{r^2(\xi+\xi_0)}{2r(\xi+\xi_0))+4}\right)^2, q\neq 0, r\neq 0, r^2-4q=0, \quad (65)$$

where, $\xi = k\left\{x - \frac{k^4(-8qr^2+r^4+16q^2)}{\Gamma(1+\alpha)}t^\alpha\right\}$.

**For Type 2:**

$$U_{1_5}(x,t) = -6k^2 pq - 6k^2 pq\tan^2(\sqrt{pq}(\xi+\xi_0)), p>0, q>0, \quad (66)$$

$$U_{1_6}(x,t) = -6k^2 pq - 6k^2 pq\cot^2(\sqrt{pq}(\xi+\xi_0)), p>0, q>0, \quad (67)$$

$$U_{1_7}(x,t) = -6k^2 pq + 6k^2 pq\tanh^2(\sqrt{-pq}(\xi+\xi_0)), p<0, q>0, \quad (68)$$

$$U_{1_8}(x,t) = -6k^2 pq + 6k^2 pq\coth^2(\sqrt{-pq}(\xi+\xi_0)), p<0, q>0, \quad (69)$$

where, $\xi = k\left\{x - \frac{k^4(16p^2q^2)}{\Gamma(1+\alpha)}t^\alpha\right\}$.

The other obtained solutions of Eq. (55) according to **Set 2** are ignored for convenience.

## 4. Discussion

The main advantage of the proposed method is that the method offers more general and huge amount of new exact traveling wave solutions with some free parameters. All the solutions obtained by the exp(-Φ(ξ))-expansion method are obtained by the applied method, and in addition one can be obtained some new explicit and exact solutions. The exact solutions have its extensive potentiality to interpret the inner structures of the natural phenomena arises in mathematical physics, chemistry and biology or any natural varied instances.

Recently, Bekir and Guner [37] have applied the ($G'/G$)-expansion method in finding the traveling wave solutions to the considered the space-time fractional Burgers equation, the space-time fractional KdV-Burgers equation and the space-time fractional coupled Burgers' equations. Gepreel and Omran [15] have employed the improved ($G'/G$)-expansion method in investigating the exact solutions to the time- and space-fractional derivative foam drainage equation. The solutions $u_{1_5}(x,t)$ and $u_{1_9}(x,t)$ of space-time fractional Burgers equation are equivalent to the solutions $u_{1,2}(x,t)$ and $U_{5,6}(\xi)$ of Bekir and Guner [37] respectively. The solutions $u_3(x,t)$, $v_3(x,t)$ and $u_1(x,t), v_1(x,t)$ of Bekir and Guner [37] to the space-time fractional coupled Burgers equation are equivalent to our obtain solutions $u_9(x,t)$, $v_9(x,t)$ and $u_7(x,t), v_7(x,t)$ respectively. Finally, the solutions $V_7(x,t)$ and $V_8(x,t)$ to the space-time fractional foam drainage equation are also coincide with the solutions in equations (24) and (23) of Gepreel and Omran [15], respectively. Very recently, Ray and Sahoo [38] have applied a novel analytical method to construct the exact solutions of time fractional fifth-order Sawada-Kotera equation. We have observed that the solution $U_{1_7}(x,t)$ is equivalent of Ray and Sahoo solution $\psi_{11}(x,t)$. Beside this, the solutions are obtained in this article via the generalized exp(-Φ(ξ))-expansion method with the auxiliary ODE (7), while the ($G'/G$)-expansion, the improved ($G'/G$)-expansion and others methods have performed with others. In this paper, the other obtained traveling waves solutions are new and have not been found in the previous literature. It noteworthy to point out that some of the obtained solutions are good in agreement with already published results, if the parameters taken particular values which authenticate to the obtained solutions. Therefore, it can be decided that the proposed method is powerful mathematical tool for solving nonlinear FPDEs. We are working on another article in which the other nonlinear FPDEs are considered for showing its effectiveness.

## 5. Conclusion

In this paper, the generalized exp(-Φ(ξ))-expansion method along with the Jumarie's modified Riemann-Liouville derivatives have proposed for solving FPDEs and applied it to seek exact traveling wave solutions to the considered equations. The obtained solutions to the considered equations may be useful for the further analysis, such as stability analysis and compare with numerical solvers arises in various areas of applied mathematics and mathematical Physics. Finally, it is worth mentioning that the implementation of this proposed method is very simple and straightforward, and it can also be applied to many other nonlinear fractional evolution equations. In future, the plan is to study the numerical simulations for these equations along with other methods. Those results will be reported in future publications.